\documentclass{article}
\usepackage{amsmath,amssymb,amsthm}
\newcommand{\set}[1]{\{#1\}}

\newcommand{\abs}[1]{\left|#1\right|}
\newcommand{\inv}{^{-1}}
\newcommand{\onto}{\twoheadrightarrow}
\newcommand{\func}[4][\to]{#2\colon#3#1#4}
\newcommand{\Z}{\mathbb{Z}}
\DeclareMathOperator{\im}{im}
\newtheorem{prop}{Proposition}
\theoremstyle{remark}
\newtheorem*{notat}{Notation}

\begin{document}
\title{Arbitrarily Large Finite Quotients Imply No Uniform Bound on Depth}
\author{Andrew~D. Warshall
\thanks{The author acknowledges support under the VIGRE program}
\\Yale University\\Department of Mathematics\\10 Hillhouse Avenue
\\P.O. Box 208283\\New Haven, CT 06520-8283\\USA}
\date{8 February 2006}
\maketitle

\begin{abstract}
We show that any group with arbitrarily large finite quotients admits
generating sets with respect to which it has arbitrarily large finite
dead-end depth. This extends a joint result with Riley and partially
answers a question asked there.
\end{abstract}

Let $G$ be any finitely generated group and $A$ a finite generating
set for $G$. Then we define the \emph{depth} (or more verbosely the
\emph{dead-end depth}) of an element $g\in G$ with respect to $A$ to
be the distance (with respect to the word metric induced by $A$) from
$g$ to the nearest (to $g$) group element at greater distance than $g$
from the identity. Formally, the depth of $g$ is the distance (in the
Cayley graph of $G$ with respect to $A$) between $g$ and the
complement of the radius-$d_A(1,g)$ closed ball about the identity,
where $d_A$ refers to distance with respect to $A$. We define the
depth of $G$ with respect to $A$ to be the (possibly infinite)
supremum of all the depths of all the elements of $G$ with respect to
$A$.

In general, the depth of a group is not the same for all finite
generating sets. In fact, in joint work with Riley \cite{RW}, we
constructed a (finitely \emph{presented}) group which has depth $2$
with respect to one generating set and depth $\infty$ with respect to
another. For hyperbolic groups, in contrast, Bogopol'sk\u{i}i \cite{B}
has shown that the depth is finite for every generating set. This
result suggests the question whether depth is bounded uniformly over
all finite generating sets. In \cite{RW}, we answered this question in
the negative for $\Z$, showing that, for all finite $n$, there exists
a generating set for $\Z$ with respect to which it has depth at least
$n$. We here answer it similarly for a broad class of groups,
including in particular all indicable groups. Specifically, we prove
the following:

\begin{prop}\label{bigdepth}
If $G$ is a finitely generated group with arbitrarily large finite
quotients and $n$ is a positive integer then $G$ has a generating set
$A$ such that the depth of $G$ with respect to $A$ is $\ge n$.
\end{prop}

This result includes the case where $G$ is indicable, for given
$G\onto\Z$ we have maps $G\onto C_n$ for arbitrarily large $n$.

Let $G$ be a finitely generated group with arbitrarily large finite
quotients and $S$ a finite generating set for $G$. We set $a=\abs{S}$.

\begin{prop}\label{biglength}
For every integer $n'$ there exists $n\ge n'$, $H_n$ and
$\func[\onto]{\pi_n}{G}{H_n}$ where $H_n$ is a finite group which has
an element $h_n\in H_n$ of length $n$ in $(\pi_n(S))^{\pm1}$ but no
element of greater length.
\end{prop}

\begin{proof}
Let $G_m$ be a quotient of $G$, $\abs{G_m}=m$, with $\alpha_m$ being
the quotient map $G\onto G_m$. If $G$ is generated by $S$ then $G_m$
is generated by $\alpha_m(S)$, which consists of $a$ elements, for all
$m$. So there are at most $(2a+1)^n$ words of length $\le n$ with
respect to that generating set, since $(2a+1)^n$ is the number of
words of length exactly $n$ in $\alpha_m(S)^{\pm1}\cup\set{1}$. In
particular, there are at most $(2a+1)^n$ group elements within $n$ of
the identity with respect to that generating set. Let $n$ be the
maximal distance from the identity of any element of $G_m$; this must
exist since $G_m$ is finite. Then $(2a+1)^n\ge m$, since all group
elements are within $m$ of the identity. So if we set $m=(2a+1)^{n'}$
then we have $n\ge n'$. Thus $G_{(2a+1)^{n'}}$ and
$\alpha_{(2a+1)^{n'}}$ will serve as $H_n$ and $\pi_n$ for some $n\ge
n'$.
\end{proof}

\begin{notat}
We denote $\pi_n(S)$ by $T_n$.
\end{notat}

\begin{proof}[Proof of Proposition~\ref{bigdepth}]
We let $N$ and $n'$ be positive integers and consider $n\ge n'$ and
the map $\func[\onto]{\pi_n}{G}{H_n}$ given by the proposition. Let
$A\subset G$ be $B_{S,N}(1)\cap\pi_n\inv(T_n)$. I claim that $G$ has
depth at least $n\ge n'$ with respect to $A$. Clearly, $A$ is a
generating set so long as $N\ge1$, since then $S\subseteq A$ since
$\pi_n(S)=T_n$. Since $h_n$ is at distance $n$ from the identity with
respect to $T_n$, any $g_n\in\pi_n\inv(h_n)$ must be at distance at
least $n$ from the identity with respect to $\pi_n\inv{T_n}$, hence
\emph{a fortiori} with respect to $A$.

But I claim that, for every $h\in H_n$, there is an element of
$\pi_n\inv(h)$ within distance $n$ of the identity with respect to
$S$. To see this, for any $h\in H_n$, take a minimal-length word
$\pi_n(s_1)\pi_n(s_2)\dots\pi_n(s_m)$ representing $h$, $s_i\in S$,
$m\le n$. Then $a=s_1s_2\dots s_m\in\pi_n\inv(h)$ is clearly at
distance at most $m\le n$ from the identity. In particular, setting
$h=h_n$, there exists $g_n\in\pi_n\inv(h_n)$ such that
$d_A(1,g_n)=n$. I claim that $g_n$ has depth at least $n$.

Let $g\in B_{A,d}(g_n)$, where $d<n$. Then $g$, being the product of
$g_n$, of length $n$, with up to $d$ words of length $N$, is, by the
triangle inequality, expressible as a word of length at most $n+dN$ in
$S$. Let $k=\abs{\pi_n(g)}$ and divide a geodesic word in $S$
representing $G$ into $k$ roughly equal pieces. Thus $g=u_1u_2\dots
u_k$, where the $u_i$ are $\in G$ and their lengths in $S$ are as
equal as possible, so that for each $u_i$ $\abs{u_i}<(n+dN)/k+1$. Let
also $t_1t_2\dots t_k$ be a minimal-length word in $T_n^{\pm1}$
representing $\pi_n(g)$. For $x\in H_n$, we denote by $\phi(x)$ a
minimal-length representative of $\pi_n\inv(x)$, so by the preceding
paragraph we always have $\abs{\phi(x)}\le n$. Then we also have
\begin{multline*}
g=[u_1\phi(\pi_n(u_1)\inv t_1)][\phi(\pi_n(u_1)\inv t_1)\inv
u_2\phi(\pi_n(u_1u_2)\inv t_1t_2)]\\ [\phi(\pi_n(u_1u_2)\inv
t_1t_2)\inv u_3\phi(\pi_n(u_1u_2u_3)\inv t_1t_2t_3)]\dots\\
[\phi(\pi_n(u_1\dots u_{k-1})\inv t_1\dots t_{k-1})\inv u_k].
\end{multline*}

We find that, if the $v_i$ are taken to be the factors as indicated by
the brackets, $g$ is now expressed as the product $v_1v_2\dots v_k$. I
claim that all the $v_i\in(\pi_n\inv(T_n))^{\pm1}$. If $i<k$, this is
clear since then
\[
\pi_n(v_i)=t_{i-1}\inv t_{i-2}\inv\dots t_1\inv\pi_n(u_1u_2\dots
u_{i-1})\pi_n(u_i)\pi_n(u_1u_2\dots u_i)\inv t_1t_2\dots t_i=t_i,
\]
while if $i=k$ then
\begin{multline*}
\pi_n(v_i)=\pi_n(v_k)=t_{k-1}\inv t_{k-2}\inv\dots
t_1\inv\pi_n(u_1u_2\dots u_{k-1})\pi_n(u_k)\\=t_{k-1}\inv
t_{k-2}\inv\dots t_1\inv\pi_n(g)\\=t_{k-1}\inv t_{k-2}\inv\dots
t_1\inv t_1t_2\dots t_k=t_k.
\end{multline*}
Furthermore, every $v_i$ is obtained from its corresponding $u_i$ by
concatenating it with up to two words each of which, being in
$\im{\phi}$, has length at most $n$, so
$\abs{v_i}\le\abs{u_i}+2n$. Thus all the $v_i$ will be $\in A^{\pm1}$
so long as all the $\abs{u_i}+2n\le N$, which will occur if
$(n+dN)/k+1+2n\le N$, for which (since $k\ge n-d$ by the triangle
inequality) it suffices that $(n+dN)/(n-d)+2n+1\le N$, which will be
true so long as $n>2d$ and $N\ge(2n^2+2nd+2n-d)/(n-2d)$.
\end{proof}

\end{document}